\documentclass[a4paper,11pt]{article}
\usepackage{amsfonts,amssymb,amsmath}

\usepackage[english]{babel}
 \makeatletter
 \@addtoreset{equation}{subsection}
 \makeatother

 \makeatletter
 \@addtoreset{enunciato}{section}
 \makeatother

 \newcounter{enunciato}[subsection]

 \newtheorem{ittheorem}{Theorem}
 \newtheorem{itlemma}{Lemma}
 \newtheorem{itproposition}{Proposition}
 \newtheorem{itcorollary}{Corollary}
 \newtheorem{itdefinition}{Definition}
 \newtheorem{itremark}{Remark}
 \newtheorem{itclaim}{Claim}
 \newtheorem{itfact}{Fact}
 \newtheorem{itconjecture}{Conjecture}

 \newenvironment{theorem}{\addtocounter{enunciato}{1}
 \begin{ittheorem}}{\end{ittheorem}}

 \newenvironment{lemma}{\addtocounter{enunciato}{1}
 \begin{itlemma}}{\end{itlemma}}

 \newenvironment{proposition}{\addtocounter{enunciato}{1}
 \begin{itproposition}}{\end{itproposition}}

  \newenvironment{corollary}{\addtocounter{enunciato}{1}
 \begin{itcorollary}}{\end{itcorollary}}

 \newenvironment{definition}{\addtocounter{enunciato}{1}
 \begin{itdefinition}}{\end{itdefinition}}

 \newenvironment{remark}{\addtocounter{enunciato}{1}
 \begin{itremark}}{\end{itremark}}

 \newenvironment{claim}{\addtocounter{enunciato}{1}
 \begin{itclaim}}{\end{itclaim}}

 \newenvironment{fact}{\addtocounter{enunciato}{1}
 \begin{itfact}}{\end{itfact}}

 \newenvironment{conjecture}{\addtocounter{enunciato}{1}
 \begin{itconjecture}}{\end{itconjecture}}

 \newcommand{\be}[1]{\begin{equation}\label{#1}}
 \newcommand{\ee}{\end{equation}}

 \newcommand{\bl}[1]{\begin{lemma}\label{#1}}
 \newcommand{\el}{\end{lemma}}

 \newcommand{\br}[1]{\begin{remark}\label{#1}}
 \newcommand{\er}{\end{remark}}

 \newcommand{\bt}[1]{\begin{theorem}\label{#1}}
 \newcommand{\et}{\end{theorem}}

 \newcommand{\bd}[1]{\begin{definition}\label{#1}}
 \newcommand{\ed}{\end{definition}}

 \newcommand{\bcl}[1]{\begin{claim}\label{#1}}
 \newcommand{\ecl}{\end{claim}}

 \newcommand{\bfact}[1]{\begin{fact}\label{#1}}
 \newcommand{\efact}{\end{fact}}

 \newcommand{\bp}[1]{\begin{proposition}\label{#1}}
 \newcommand{\ep}{\end{proposition}}

 \newcommand{\bc}[1]{\begin{corollary}\label{#1}}
 \newcommand{\ec}{\end{corollary}}

 \newcommand{\bcj}[1]{\begin{conjecture}\label{#1}}
 \newcommand{\ecj}{\end{conjecture}}

 \newcommand{\bpr}{\begin{proof}}
 \newcommand{\epr}{\end{proof}}

\newcommand{\bprp}[1]{\begin{proofofp}{\it\ref{#1}}.\,\,}
 \newcommand{\eprp}{\end{proofofp}}

 \newcommand{\bprl}[1]{\begin{proofofl}{\it\ref{#1}}.\,\,}
 \newcommand{\eprl}{\end{proofofl}}

 \newcommand{\bi}{\begin{itemize}}
 \newcommand{\ei}{\end{itemize}}

 \newcommand{\ben}{\begin{enumerate}}
 \newcommand{\een}{\end{enumerate}}

 \newenvironment{proof}{\noindent {\em Proof}.\,\,}{\hspace*{\fill}$\halmos$\medskip}
 \newenvironment{proofofp}{\noindent {\em Proof of Proposition\,\,}}{\hspace*{\fill}$\halmos$\medskip}
 \newenvironment{proofofl}{\noindent {\em Proof of Lemma\,\,}}{\hspace*{\fill}$\halmos$\medskip}
 \newcommand{\halmos}{\rule{1ex}{1.4ex}}

 \newcommand{\one}{{\mathchoice {1\mskip-4mu\mathrm l}
         {1\mskip-4mu\mathrm l}
         {1\mskip-4.5mu\mathrm l}
         {1\mskip-5mu\mathrm l}}}

\def \Z {{\mathbb Z}}

\def \ra {\rightarrow}
\def \ba {\begin{array}}
\def \ea {\end{array}}

\def \ua {\uparrow}

\def \P {{\mathbb P}}
\def \E {{\mathbb E}}

\def \Var {{\rm Var}}
\def \Cov {{\rm Cov}}

\def \c {{\rm c}}

\def \cF {{\mathcal F}}

\def\one{\rlap{\mbox{\small\rm 1}}\kern.15em 1}

\def\embf#1{\emph{\bf #1}}

\begin{document}
\title{Large deviations for voter model occupation times in two dimensions\\
\bigskip
Grandes d\'eviations pour le temps d'occupation du mod\`ele du votant bidimensionnel}

\author{\renewcommand{\thefootnote}{\arabic{footnote}}
G.\ Maillard
\footnotemark[1]
\\
\renewcommand{\thefootnote}{\arabic{footnote}}
T.\ Mountford
\footnotemark[1]
}

\footnotetext[1]{
Institut de Math\'ematiques, \'Ecole Polytechnique F\'ed\'erale
de Lausanne, Station 8, CH-1015 Lausanne, Switzerland,
{\sl gregory.maillard@epfl.ch},
{\sl thomas.mountford@epfl.ch}
}
\date{22th April 2008}
\maketitle

\begin{abstract}
We study the decay rate of large deviation probabilities of occupation times,
up to time $t$, for the voter model $\eta\colon\Z^2\times[0,\infty)\ra\{0,1\}$
with simple random walk transition kernel, starting from a Bernoulli product
distribution with density $\rho\in(0,1)$. In \cite{bramcoxgri88}, Bramson, Cox
and Griffeath showed that the decay rate order lies in $[\log(t),\log^2(t)]$.

In this paper, we establish the true decay rates depending on the level. We show
that the decay rates are $\log^2(t)$ when the deviation from $\rho$ is maximal
(i.e., $\eta\equiv 0$ or $1$), and $\log(t)$ in all other situations. This answers
some conjectures in \cite{bramcoxgri88} and confirms nonrigorous analysis carried
out in \cite{benfrakra96}, \cite{dorgod98} and \cite{howgod98}.

\begin{center}
{\bf R\'esum\'e}
\end{center}

\vspace{-0.1cm}
On \'etudie le taux de d\'ecroissance des probabilit\'es de grandes d\'eviations
des temps d'occupation, jusqu'\`a l'instant $t$, du mod\`ele du votant
$\eta\colon\Z^2\times[0,\infty)\ra\{0,1\}$ ayant le noyau de transition
d'une marche al\'eatoire simple et partant d'une distribution
produit de Bernoulli de param\`etre $\rho\in(0,1)$. Dans \cite{bramcoxgri88},
Bramson, Cox et Griffeath ont montr\'e que l'ordre du taux de d\'ecroissance
se situe dans $[\log(t),\log^2(t)]$.

Dans cet article, nous \'etablissons les taux de d\'ecroissance exacts d\'ependant
du niveau. On prouve que les taux de d\'ecroissance sont $\log^2(t)$ lorsque la
d\'eviation de $\rho$ est maximale (i.e., $\eta\equiv 0$ ou $1$), et $\log(t)$
dans toutes les autres situations. Ceci r\'epond \`a une conjecture de
\cite{bramcoxgri88} et confirme l'analyse non rigoureuse effectu\'ee dans
\cite{benfrakra96}, \cite{dorgod98} et \cite{howgod98}.

\vskip 1truecm
\noindent
{\it MSC} 2000. Primary 60F10, 60K35; Secondary 60J25.\\
{\it Key words and phrases.} Voter model, large deviations.\\
{\it Acknowledgment.} SNSF: Grant $\#\, 200021-107425/1$.
The second author benefited from a visit to Centro Ennio De
Giorgi in Pisa and a conversation there with Claudio Landim.
We are grateful for comments made by Maury Bramson.
\end{abstract}

\newpage

%%%%%%%%%%%%%% SECTION 1 %%%%%%%%%%%%%%%%%%%%%%%%%%%%%%%%%%%%%%%%%%%%%%%%%%%%%%%%%%%%%%%%%%%%%%%%%%

\section{Introduction and main results}
\label{S1}

%%%%%%%%%%%%%%%%%% SUBSECTION 1.1 %%%%%%%%%%%%%%%%%%%%%%%%%%%%%%%%%%%%%%%%%%%%%%%%%%%%%%%%%%%%%%%%%

\subsection{The Voter Model}
\label{S1.1}

We consider the simple voter model in $\Z ^{2}$ corresponding to the simple random walk.
In general dimensions this voter model is a Markov process on $\{0, 1\}^{\Z^d}$ with
operator
\be{intr1}
\Omega f(\eta)
=\frac{1}{2d}\sum_{x\in\Z^{d}}\sum_{y\sim_n x}\Big(f(\eta^{x,y})- f(\eta)\Big),
\ee
where $x\sim_n y$ means $x$ and $y$ are nearest neighbours on the $\Z^{d}$ lattice and
$\eta^{x,y}$ is the configuration
\be{intr3}
\begin{cases}
\eta^{x,y} (z)= \eta (z)& \textrm{for } z\not=x,\\
\eta^{x, y} (x)= \eta(y).
\end{cases}
\ee

This process was introduced independently by Clifford and Sudbury \cite{clisud73} and by Holley
and Liggett \cite{hollig75}.  There the basic results concerning equilibria were shown: for
recurrent random walks (i.e.\ $d \leq 2$) the only extremal equilibria are $\delta_{0}$ and
$\delta_{1}$ whereas for transient random walks there exists for each $\rho\in [0,1]$ an
extremal, translation invariant ergodic equilibrium of density $\rho$, $\mu_{\rho}$ (and these
are the totality of extremal equilibria). In the transient case the measures $\mu_{\rho}$ are
the limits for distributions of the process begun with initial measure $\nu_{\rho}$ for which
$(\eta(x):x \in \Z^{d})$ are i.i.d.\ Bernoulli $(\rho)$ random variables. Details for this
and much more can be found in Liggett \cite{lig85}.

In this note our analysis will rely heavily on the duality of the voter model with coalescing
random walks (as exploited in \cite{bramcoxgri88} and \cite{cox88}--\cite{coxgri86}): given
distinct space time points ${(x_{i},t_{i})}_{i=1}^{r}$ in $\Z^{d}\times[0,\infty)$, the joint
distribution of $(\eta_{t_{i}}(x_{i}))_{i =1}^{r}$ can be determined via coalescing random
walks $(\chi_{t}^{i}:t\geq 0)$ defined as follows: (suppose without loss of generality that
$0\leq t_{1}\leq t_{2} \dots\leq t_{r}$) $\chi_{t}^{i}= x_{i}$ for $0\leq t\leq t_{r}-t_{i}$,
thereafter $\chi^{i}$ evolves as a simple random walk. If for $i<j$, $s\geq t_{r}-t_{i}$,
$\chi_{s}^{i}=\chi_{s}^{j}$, then $\chi_{s'}^{i}=\chi_{s'}^{j}$ for all $s^{'}\geq s$. That is
the random walks are coalescing. Otherwise the random walks evolve independently. The joint law
of $(\eta_{t_{1}}(x_{1}),\eta_{t_{2}}(x_{2}),\ldots,\eta_{t_{r}}(x_{r}))$ is that of
$(\eta_{0}(\chi^{1}_{t_{r}}),\eta_{0}(\chi^{2}_{t_{r}}),\ldots,\eta_{0}(\chi^{r}_{t_{r}}))$.
The clear exposition of the Harris construction of the voter model found in Durrett \cite{D1} is
here recommended.

We will in this article be concerned with the behaviour, for $t$ large, of
\be{intr5}
\frac{T_{t}}{t}=\frac{1}{t}\int_{0}^{t}\eta_{s}(0)\, ds
\ee
for $(\eta_{s}:s\geq 0)$ a voter model begun with initial measure $\nu_{\rho}$, $\rho\in(0,1)$.
In the transient regime, the behaviour is equivalent to that for a voter model begun with initial
distribution $\mu_{\rho}$. This problem was discussed in a series of papers by Cox and Griffeath
\cite{coxgri83} and \cite{coxgri86} and Bramson Cox and Griffeath \cite{bramcoxgri88} .
It follows from the duality description also, as noted in these articles, that $T_{t}$ may be
understood as follows.

A Harris system for the voter model $(\eta_t: t \geq 0)$ is a collection of independent rate
$\frac{1}{2d}$ Poisson processes $N^{x,y}$ for every ordered pair $x,y$ with $y \sim_n x$.
From this system $ \eta_.$ evolves by stipulating that for $x \in \Z^d , \eta_t(x) $ changes
value (or flips) only at times $t$ in $N^{x,y}$ for some $y \sim_n x$.  At such a time $t$ we
put $\eta_t(x) = \eta_t(y)$.  If for $t \in N^{x,y}, \eta_{t-}(x) = \eta_t(y) $ then there
is no change in value for $ \eta_.(x) $ at time $t$.  Given this system we can define for each
$x \in \Z^d$ and $t \geq 0$ dual simple random walks, $(Z^{x,t}_s: 0 \leq s \leq t)$ with
$Z^{x,t}_0=x$, as follows:
\be{har3}
Z^{x,t}_s \ne \ Z^{x,t}_{s-} \iff \ t-s \in N^{Z^{x,t}_{s-},w}
\ee
for some $w \sim_n Z^{x,t}_{s-}$.  In which case $Z^{x,t}_.$ jumps from $Z^{x,t}_{s-}$
to $w$ at time $s$.

The importance of these random walks lies in the following properties\\
1) $\eta_t(x) = \eta_0(Z^{x,t}_{t})$\\
and \\
2) the random walks $\{Z^{x,t}_{.}\}_{x \in \Z^d}$ are independent until they meet.

Furthermore it may be seen that if $0 \leq t \leq t ^\prime $ then the random walks
$Z^{x,t}_{.}$ and $Z^{x,t^\prime }_{.}$ are coalescing in the sense that \\
3) if for some $s \in [0,t], \ Z^{x,t}_{s} = Z^{x,t ^\prime }_{t^\prime - t+s}$, then
$Z^{x,t}_{u} = Z^{x,t ^\prime }_{t^\prime - t+u}$ for all $s \leq u \leq t$.

Thus we have
 a system of coalescing random walks $(\chi_{v}^{s},v\in[0,s]) = Z^{0,s}_{v}$ on
$\Z^{d}$, so that $\chi_{0}^{s}=0$ and by property 3) above if for $0 \leq  v \leq s \leq s^\prime,
\chi^{s}_{v} = \chi^{s ^\prime }_{s^\prime - s+v}$, then $ \chi^{s}_{u} =
\chi^{s ^\prime }_{s^\prime - s+u}
$ for all $v \leq u \leq s$.  We call the collection of random walks $\{\chi^s_. \}_{s \geq 0}$ the
coalescing random walks associated with $\eta_.(0)$.

Let $O_{x} = \lambda_t(\{s \in [0, t] : \chi_{s}^{s} = x\})$ with $\lambda_t$ the Lebesgue measure
on $[0,t]$, then
\be{intr7}
T_{t}
=\sum_{x\in\{\chi_{s}^{s}:s\in[0, t]\}} O_{x}\, \eta_{0} (x).
\ee
For $\eta_{0}$ distributed as product measure $\nu_{\rho}$ the duality representation immediately
yields
\be{intr9}
\begin{aligned}
\Var \bigg(\frac{T_{t}}{t}\bigg)
&= \frac{1}{t^{2}} \int_{0}^{t}\int_{0}^{t} \Cov \big(\eta_{s}(0), \eta_{s'}(0)\big)\, ds'\, ds\\
&= \frac{1}{t^{2}}\int_{0}^{t}\int_{0}^{t} P\Big(\chi_{s}^{s}=\chi_{s'}^{s'}\Big) \rho(1-\rho)\,
ds'\,ds.
\end{aligned}
\ee
For $s>s'$,  $P(\chi_{s}^{s}=\chi_{s'}^{s'})$ is easily seen to be the probability that a random walk
issuing from the origin hits the origin during the interval $[s-s',s+s']$. If one chooses $s,s'$
uniformly on $[0,t]$ this probability is easily seen to tend to zero as $t\to\infty$ for transient
random walks. However for recurrent random walks it may tend to zero as $t\to\infty$ (for $d = 2$) or
it may tend to a non zero limit $(d=1)$. From this we obtain: for $\eta_{0}$ distributed by
$\nu_{\rho}$,
\be{intr11}
\frac{T_{t}}{t} \longrightarrow \rho \quad \textrm{in probability if only if}\quad d \geq 2.
\ee
In fact, the convergence in (\ref{intr11}) holds a.s. (see Cox and Griffeath \cite{coxgri83}).

%%%%%%%%%%%%%%%%%% SUBSECTION 1.2 %%%%%%%%%%%%%%%%%%%%%%%%%%%%%%%%%%%%%%%%%%%%%%%%%%%%%%%%%%%%%%%%%

\subsection{Asymptotic behavior of occupation time}
\label{S1.2}

 Bramson, Cox and Griffeath \cite{bramcoxgri88} obtained large
deviation bounds: for each $\alpha\in(\rho,1]$ there exist positive finite constants
$C_1=C_1(d)$, $C_2=C_2(d,\alpha)$ such that, for $t$ sufficiently large,
\be{ratevm}
\begin{cases}
e^{-C_1 \log^2(t)} \leq \P_{\nu_\rho}\big(T_t \geq \alpha t\big)
\leq e^{-C_2 \log(t)}
&\text{if } d=2,\\[0.3cm]
e^{-C_1 b_t} \leq \P_{\nu_\rho}\big(T_t \geq \alpha t\big) \leq e^{-C_2 b_t}
&\text{if } d\geq 3,\\
\end{cases}
\ee
with
\be{ratevm2}
b_t =
\begin{cases}
\sqrt{t}
&\text{if } d=3,\\
\frac{t}{\log t}
&\text{if } d=4,\\
t
&\text{if } d \geq 5.\\
\end{cases}
\ee
By symmetry arguments, the large deviation regime is the same for the
deviations $T_t/t\leq \alpha$ with $\alpha\in[0,\rho)$.

%%%%%%%%%%%%% SECTION 1.3 %%%%%%%%%%%%%%%%%%%%%%%%%%%%%%%%%%%%%%%%%%%%%%%%%%%%%%%%%%%%%%%%%%%%%%%%%

\subsection{Results}
\label{S1.3}
Given the bounds of \cite{bramcoxgri88} cited in the previous section, in so far as the
exponential order of large deviations is concerned, the only outstanding case is the
two-dimensional one. Throughout the rest of the paper, we assume that $d=2$.  The
following two results constitute a full resolution of the question of exponential order
for the large deviations of $T_t$.
\bt{Th1}
There exist positive finite constants $C_1, C_2$ such that, for t sufficiently large,
\be{Th1bounds}
e^{-C_1 \log^2(t)}
\leq \P_{\nu_\rho}\big(T_t=t\big)
\leq e^{-C_2 \log^2(t)}.
\ee
\et

\bt{Th2}
For each $\alpha\in(\rho,1)$, there exist positive finite constants $C_1=C_1(\alpha)$,
$C_2=C_2(\alpha)$ such that, for t sufficiently large,
\be{Th2bounds}
e^{-C_1 \log(t)}
\leq \P_{\nu_\rho}\big(T_t\geq\alpha t\big)
\leq e^{-C_2 \log(t)}.
\ee
\et
By (\ref{ratevm}), it only remains to prove the upper bound in Theorem \ref{Th1} and
lower bound in Theorem \ref{Th2}. If $g(t)$ and $h(t)$ are real functions, we write
$g(t)\asymp h(t)$ as $t\to\infty$ when
\be{sim}
0<\liminf_{t\to\infty}g(t)/h(t)\leq\limsup_{t\to\infty}g(t)/h(t)<\infty.
\ee

%%%%%%%%%%%%% SECTION 2 %%%%%%%%%%%%%%%%%%%%%%%%%%%%%%%%%%%%%%%%%%%%%%%%%%%%%%%%%%%%%%%%%%%%%%%%%%%

\section{Discussion}
\label{S2}

The study of $T_t$ was initiated by Cox and Griffeath \cite{coxgri83} who noted that the
question of its large deviations belonged naturally with related issues arising in the
Ising and percolation models, but that in contrast (and due to the tractable duality) with
these, progress in identifying the effect at low dimensions was possible. Nevertheless
questions remain.

The behavior of $T_t$ in low dimensions has motivated studies in the Physics community.
Due to the recurrence of simple random walks, as $t\to\infty$, the simple voter model forms
larger and larger clusters when $d\leq 2$ (a more detailed analysis of clustering can be
found in \cite{coxgri86}). Therefore, a consensus of opinion is approached as $t\to\infty$.
In words, that means that the system coarsens.  A natural question is to study, for such a
corsening system, the asymptotic behavior of the persistence probability $\P\big(T_t=t\big)$,
i.e., the probability that a given site will never change its state as time goes to infinity.
To be in accordance with the physicist terminology, consider the voter model
$\zeta\colon\Z^d\times[0,\infty)\ra\{-1,1\}$ (as a spin system) with opinions $-1$ and $+1$.
Define the \emph{mean magnetization} at time $t$ by
\be{meanmagn}
M(t)=\frac1t\int_0^t \zeta(0,s)\, ds,
\quad M(t)\in[-1,1].
\ee
In the case considered the initial distribution was symmetric w.r.t.\ $-1$ and $1$ and
so $\E(M(t))=0$.
Then for all $x>0$, the distribution of the mean magnetization, $P(t,x)=\P(M(t)\geq x)$,
and $R(t,x)=P(M(s)\geq x,\,\,\forall\,s\leq t)$ represent the \emph{deviation} of $M(t)$
from its mean and the probability of \emph{persistent large deviations}, respectively.
Then, assuming that $\zeta(0,0)=1$,
\be{persistprob}
P(t,1)=R(t,1)=\P\big(\zeta(0,s) = 1 \,\,\forall\, 0\leq s\leq T\big)
\ee
is the so called \emph{persistent probability} and corresponds to the object of study
of Theorem \ref{Th1}. Ben-Naim, Frachebourg and Krapivsky \cite{benfrakra96} showed
convincingly via numerical methods that there exists some $C>0$ such that
\be{persistprob1}
P(t,1)\asymp e^{-C \log^2(t)},
\qquad t\gg 1.
\ee
Howard and Godr\`eche \cite{howgod98} confirm nonrigorously this result both by using
path-integral methods and Monte Carlo simulations. After a sharper analysis, Dornic and
Godr\`eche \cite{dorgod98} concluded that
\be{persistprob2}
P(t,x)\asymp e^{-I(x) \log(t)}
\quad\text{and}\quad
R(t,x)\asymp e^{-J(x) \log^2(t)}
\qquad t\gg 1
\ee
with $\lim_{x\to 1}I(x)=\infty$ and $\lim_{x\to 1}J(x)=C$ for some constant $C>0$.
This is in accordance with Theorems \ref{Th1}--\ref{Th2}.

%%%%%%%%%%%%% SECTION 3 %%%%%%%%%%%%%%%%%%%%%%%%%%%%%%%%%%%%%%%%%%%%%%%%%%%%%%%%%%%%%%%%%%%%%%%%%%%

\section{Proofs}
\label{S3}

%%%%%%%%%%%%% SUBSECTION 3.1 %%%%%%%%%%%%%%%%%%%%%%%%%%%%%%%%%%%%%%%%%%%%%%%%%%%%%%%%%%%%%%%%%%%%%%

\subsection{Proof of Theorem \ref{Th1}}
\label{S3.1}

Let $\chi=(\chi^t)_{t\geq 0}=(\chi_s^t, s\in [0,t])_{t\geq 0}$ be the coalescing random walks
associated with $\eta_.(0)$ for a voter model $(\eta_t: t \geq 0)$.  Denote by $P$ and $E$,
respectively, probability and expectation associated with $\chi$. The dual relationship between
voter model and coalescing random walks lead to the following lemma (see Bramson, Cox and
Griffeath \cite{bramcoxgri88}, Section 1 for details).
\bl{chilem}
For all $t>0$
\be{chi3}
\P_{\nu_\rho}\big(T_t=t\big) =E\Big(\rho^{\#\chi^t}\Big),
\ee
where $\#\chi^t$ denote the number of distinct sites in the collection $\{ \chi^s_s : 0 \leq s \leq t \}$.
\el
Then, the proof of Theorem \ref{Th1} reduces to the following proposition.
\bp{Prop1}
There exist $K_1,K_2>0$ so that
\be{chi5}
P\big(\#\chi^t\leq K_1\log^2(t) \big) \leq e^{-K_2 \log^2(t)}
\ee
for all $t>0$ sufficiently large.
\ep
Indeed, combining Lemma \ref{chilem} and Proposition \ref{Prop1}, we get
\be{chi7}
\begin{aligned}
\P_{\nu_\rho}\big(T_t=t\big)
&= E\Big(\rho^{\#\chi^t}\one\big\{\#\chi^t\leq K_1\log^2(t)\big\}\Big)
+E\Big(\rho^{\#\chi^t}\one\big\{\#\chi^t > K_1\log^2(t)\big\}\Big)\\
&\leq P\Big(\#\chi^t\leq K_1\log^2(t)\Big) + \rho^{K_1\log^2(t)}\\
&\leq  e^{-C_2 \log^2(t)},
\end{aligned}
\ee
where in the last inequality we choose $K_1$ small enough and $t$ sufficiently large.
This completes the proof of Theorem \ref{Th1}. The next section is devoted to the proof
of Proposition \ref{Prop1}.

%%%%%%%%%%%%% SUBSECTION 3.2 %%%%%%%%%%%%%%%%%%%%%%%%%%%%%%%%%%%%%%%%%%%%%%%%%%%%%%%%%%%%%%%%%%%%%%

\subsection{Proof of Proposition \ref{Prop1}}
\label{S3.2}

The overall strategy is to show that on an interval $[3t/4,t]$ with probability
of the order $1-e^{-C_1\log(t)}$ for some universal $C_1>0$, the stream of coalescing
random walks produces $C_1\log(t)$ distinct random walks which hit the annulus
$B(0,\sqrt{2t})\setminus B(0,\sqrt{t})$, where $B(0,t)=\{x\in\Z^2\colon |x|\leq t\}$
($t\geq 0$), before time $t/2$ and do not leave in dual time $[t/2,t]$. If we call this
event $A_t$, then it can be shown that $A_t,A_{t/2},A_{t/4},\ldots$ are independent,
each producing with probability $1-e^{-C_1\log(t)}$, of the order $\log(t)$ distinct
random walks. This will be enough to show Proposition \ref{Prop1}.

In order to prove Proposition \ref{Prop1}, we need a number of preparatory results
concerning ordinary and coalescing random walks.
Let $X=(X(u)\colon u\geq 0)$ be a simple random walk on $\Z^2$ with continuous time
transition probability kernel $p_u(\,\cdot\,)$. Denote by $P^x$ its probability
law starting from $x\in\Z^2$ and for all $y\in\Z^2$, $t>0$, let
\be{hitting}
\tau_y=\inf\{u>0\colon X(u)=y\}
\quad\text{and}\quad
\sigma_t=\inf\{u>0\colon |X(u)|\geq t\}.
\ee
We refer to Lawler \cite{law91} for hitting probabilities for the two dimensional
simple random walk:
\bl{2DRWbdslem}
Uniformly for $x\in\Z^2\setminus\{0\}$, $|x|\leq \sqrt{t}$,
\be{RWbd1}
P^{x}\Big(\tau_0<\sigma_{\sqrt{t}}\Big)
\asymp \frac{\log(\sqrt{t})-\log(|x|)}{\log(\sqrt{t})}
\quad\text{as}\quad
t\to\infty,
\ee
and
\be{RWbd2}
P^{x}\left(\tau_0<t\right)
\asymp \frac{\log(\sqrt{t})-\log(|x|)+1}{\log(\sqrt{t})}
\quad\text{as}\quad
t\to\infty.
\ee
\el
\bpr
The proof can be found in Lawler \cite{law91}, Proposition 1.6.7 in the case
of discrete time random walks. The transfer to continuous time is easy.
\epr

We now consider two independent simple random walks $\{X(u):u\geq 0\}$ and
$\{Y(u):u\geq s\}$, both starting from $0$ in the sense that $X(0)=0=Y(s)$.
We are interested in the probability that
\be{Adef}
\Big\{\exists\, s\leq u\leq t : X(u)=Y(u)\Big\} := A^{X,Y}(s,t).
\ee
\bl{Abdslem}
There exists positive constants $K_3,K_4$ so that for $s\in (t/\log(t), t/2)$
and $t$ large,
\be{Abds3}
K_3\frac{\log(t)-\log(s)}{\log(t)}
\leq P\Big(A^{X,Y}(s,t)\Big)
\leq K_4\frac{\log(t)-\log(s)}{\log(t)}.
\ee
\el
\bpr

We first show the lower bound $P(A^{X,Y}(s, t))$. We condition on the value of $X(s)$.
Thus, $P(A^{X,Y}(s,t)\,\,\vert\,\, X(s) = x)$ is equal to $P^{x}(\tau_{0}<2(t-s))$, since
$(Y(s+u)-X(u))_{u\geq 0}$ is a speed two random walk. Then given the constraints on $s$
we have
\be{}
P\Big(A^{X,Y}(s,t)\,\,\big\vert\,\, X(s)=x\Big)\geq P^{x}\big(\tau_{0}<t\big).
\ee
So
\be{Abds7}
\begin{aligned}
P\Big(A^{X,Y} (s, t)\Big)
&\geq\sum_{|x|\leq \frac{\sqrt{s}}{2}}P\big(X(s)=x\big)\, P^{x}(\tau_{0}<t)\\
&\geq C\sum_{|x|\leq\frac{\sqrt{s}}{2}} P\big(X(s)=x \big)\,
\frac{\log(t)-\log(s)}{\log(t)},
\end{aligned}
\ee
by Lemma \ref{2DRWbdslem}, for universal strictly positive $C$. This in turn is
\be{Abds9}
\geq CC'\, \frac{\log(t)-\log(s)}{\log(t)},
\ee
by the central limit for $Y(s)$. For the opposite inequality we obtain, arguing similarly,
that
\be{Abds11}
P\bigg(A^{X,Y}(s,t)\cap\bigg\{|Y(s)|\geq\frac{\sqrt{s}}{2}\bigg\}\bigg)\leq C''\,
\frac{\log(t)-\log(s)}{\log(t)}.
\ee
So it suffices to bound appropriately
\be{Abds13}
\begin{aligned}
P\bigg(A^{X,Y}(s,t)\cap\bigg\{|Y(s)|<\frac{\sqrt{s}}{2}\bigg\}\bigg)
&=\sum_{i=1}^{\lceil \log_2(\sqrt{s})\rceil}P\Big(A^{X,Y}(s,t)
\cap\Big\{|Y(s)|\in \big[\sqrt{s}2^{-i-1},\sqrt{s} 2^{-i}\big)\Big\}\Big)\\
&\qquad+P\big(Y(s)=0 \big).
\end{aligned}
\ee
Given the condition that $s\leq t/2$,
\be{Abds15}
\frac{\log(t)-\log(s)}{\log(t)}\geq\frac{\log(2)}{\log(t)}\gg P\big(Y(s)=0\big)
\ee
for $s\geq t/\log(t)$, so we may ignore the term $P(Y(s)=0)$.
By the local central limit theorem (see e.g. Durrett \cite{dur05}),
\be{Abds17}
P\Big(Y(s)\in \big[\sqrt{s}\, 2^{-i-1},\sqrt{s}\, 2^{-i}\big)\Big)
\leq K \, 2^{-2 i}
\ee
for universal $K$. By Lemma \ref{2DRWbdslem} and given the condition that
$s\in(t/\log(t),t/2)$,
\be{Abds19}
P\Big(A^{X,Y}(s,t)\,\,\big\vert\,\, |Y(s)|\in \big[\sqrt{s}2^{-i-1},\sqrt{s} 2^{-i}\big)\Big)
\leq \frac{1}{\log(t)}\Big(\log(t)-\log(s)+(2i+3)\log(2)+2\Big).
\ee
Combining (\ref{Abds17}--\ref{Abds19}), we get
\be{Abds21}
\begin{aligned}
&\sum_{i = 1}^{\lceil\log _2(\sqrt{s})\rceil}P\Big(A^{X,Y}(s,t)
\cap \Big\{|Y(s)|\in\big[\sqrt{s}\, 2^{-i-1},\sqrt{s}\, 2^{-i}\big)\Big\}\Big)\\
&\quad\leq  K\sum_{i = 1}^{\lceil\log _2(\sqrt{s})\rceil} 2^{-2i}\,
\frac{\log(t)- \log(s)+(2i+3)\log(2)+2}{\log(t)}\\
&\quad\leq K'\, \frac{\log(t)-\log(s)}{\log(t)},
\end{aligned}
\ee
for some $K'>0$ and we are done.
\epr

\bc{cor1}
Given $C>1$ let
\be{Rdef}
R=\bigg\lceil\frac{\log(t)}{5C}\bigg\rceil
\ee
and let $(Y^k(t) \colon t\geq t_k)$, $0\leq k\leq R$ be independent random
walks starting at
\be{meanV}
Y^k(t_k) = 0
\quad\text{with}\quad
t_k=\frac{kCt}{\log(t)}.
\ee
Then, there exists some universal (not depending on $C$) strictly positive $K_5$
so that, for all $t$ sufficiently large,
\be{meanV3}
E(V)
\leq \frac{K_5}{C}
\quad\text{with}\quad
V = \sum_{k=1}^{R}\one\Big\{A^{Y^0,Y^k}\big(t_k,t\big)\Big\}.
\ee
\ec
\br{rem1}
$E(V\,\vert\,Y^0)$ is a functional of the random walk path independent
of the random walks $Y^k$, $1\leq k\leq R$, and can and will be considered as defined for any
random walk starting at the origin, see Corollary \ref{cor2}.
\er

\bpr
By Lemma \ref{Abdslem}, for all $t$ sufficiently large
\be{meanV5}
\begin{aligned}
E(V)
\leq -K_5 \sum_{k=1}^{\big\lceil\frac{\log(t)}{5C}\big\rceil}
\frac{\log\left(\frac{kC}{\log(t)}\right)}{\log(t)}
&\leq -\frac{K_5}{C}\int_0^{\frac{\log(t)}{5C}+1}
\frac{C}{\log(t)}\log\left(\frac{Cx}{\log(t)}\right)\, dx\\
&\leq-\frac{K_5}{C}\int_0^{1} \log(x)\, dx
= \frac{K_5}C.
\end{aligned}
\ee
\epr

We now collect a few nice properties of our random walks: let $(X(u)\colon u\geq 0)$ be
a simple random walk starting at $X(0)=0$. For $t\geq 0$, recall that
$B\big(0,t\big)=\{x\in\Z^2 \colon |x|\leq t\}$.
\bl{hitballlem}
For all $t\geq 0$ and for whatever finite choice of $C\geq 1$,
\be{hitball1}
P^0\Big(X(u) \in B\big(0,t^{1/3}\big) \text{ for some }
u\in \big(t_1-1,t\big)\Big)\longrightarrow 0
\quad\text{as}\quad t\ra\infty,
\ee
for $t_1\,=\, Ct/\log(t)$.
\el
\br{rm2}
We will explain the choice of $t_1-1$ later (see Remark \ref{rm4}).
\er

\bpr
First, remark that
\be{hitball3}
P\bigg(X\big(t_1\big)\geq \frac{\sqrt{t}}{\log(t)}\bigg)\ua 1
\quad\text{as}\quad t\to\infty.
\ee
For any random process $(Z(t) \colon t\geq 0)$ on $\Z^2$ denote
\be{hitball5}
S_{t}(Z)=\inf\{s \colon |Z(s)|\leq t\}.
\ee
Therefore, in order to prove (\ref{hitball1}), it suffices to prove that
for all $|x|\geq\sqrt{t}/\log(t)$
\be{hitball7}
P^x\Big(S_{t^{1/3}}(X)<t\Big)\to 0
\quad\text{as}\quad t\to\infty.
\ee
But this follows from random walks embedding into Brownian motions and the fact that
(\ref{hitball7}) is fulfilled when a two dimensional Brownian motion is considered instead
of $X$.
\epr

The following is simply a consequence of the invariance principle.
\bl{cv2BMlem}
As $t\ra\infty$,
\be{cv2BM3}
P^0\left(|X(u)|\in \big(\sqrt{t},\sqrt{2t}\big)\,\,\, \forall\,
\frac{t}4\leq u\leq t\right)
\ra P\left(|B(u)|\in(1,\sqrt{2})\,\,\, \forall\, \frac14\leq u\leq 1\right)
=\alpha>0,
\ee
where $B$ denotes a standard two dimensional Brownian motion.
\el
We are ready to choose our constant $C$: we choose $C$ so that for $K_5$ as in
Corollary \ref{cor1} and $\alpha$ as above,
\be{Cconstant}
\frac{K_5}{C}\leq \frac{\alpha^2}{10^4}.
\ee
\bc{cor2}
For $E(V\,\vert\,X)$ as defined in Remark \ref{rem1} and $t$ sufficiently large,
the probability that the path $\{(u,X(u)) \colon 0\leq u\leq t\}$ is such that either
$$
\begin{aligned}
\text{(i)} &\quad E\big(V \,\,\vert\,\, X\big)\geq \alpha/10^2,\\
\text{or (ii)} &\quad \big|X(u)\big|\notin (\sqrt{t},\sqrt{2t}) \text{ for some } u\in(t/4,t],\\
\text{or (iii)} &\quad \big|X(u)\big|<t^{1/3} \text{ for some } u\in [t_1-1, t],
\end{aligned}
$$
is at most $1-2\alpha/3$.
\ec
\bpr
By Corollary \ref{cor1} and our choice of $C$, we have
\be{cor2-5}
P\bigg(E\big(V\,\,\vert\,\, X\big)\geq\frac{\alpha}{10^2}\bigg)
\leq \frac{10^2K_5}{\alpha C}\leq \frac{\alpha}{10^2}.
\ee
Then, combining Lemmas \ref{hitballlem}--\ref{cv2BMlem} and (\ref{cor2-5}),
we get the claim.
\epr

We consider the system of coalescing random walks $(X^i(s)\colon t_i\leq s\leq t,\,
0\leq i\leq R) = (\chi_{s-t_i}^{t-t_i}\colon t_i\leq s\leq t,\, 0\leq i\leq R)$.
We are interested in the number of distinct random walks at time $t$ which satisfy
\be{RWgood}
\left|X^i(u)\right|\in \big(\sqrt{t},\sqrt{2t}\big)
\quad\forall\, u\in\bigg[\frac{t}2, t\bigg],
\ee
where $X^i$, $0\leq i\leq R$, are coalescing random walks defined in (\ref{meanV}).
We will in turn let the random walks evolve until something ``bad'' happens. This will
mean the violation of some given conditions: Define times
\begin{itemize}
\item[(a)]
$T^{i,a}=\inf\Big\{s\geq t_{i+1} \colon \sum_{j=i+1}^{R}
P\Big( X^j(v) = X^i(v)\text{ for some }v\in\big[t_{i+1},s\big] \,\,\big\vert\,\, X^i\Big)
\geq \alpha/10^2\Big\}$;
\item[(b)] $T^{i,b}=\inf\Big\{s\geq t_{i+1}-1 \colon |X^i(s)|\leq t^{1/3}\Big\}$;
\item[(c)] $T^{i,c}=\inf\Big\{s\geq t_i+t/4 \colon |X^i(s)|\notin (\sqrt{t},\sqrt{2t})\Big\}$;
\item[(d)] $T^i=t\wedge T^{i,a}\wedge T^{i,b}\wedge T^{i,c}$,
\end{itemize}
and kill (or freeze) the random walk $X^i$ at time $T^i$.

\br{rm3}
Note that in (c), since for all $0\leq i\leq R$, $t_i\leq t/4$, $X^i$ will satisfy
(\ref{RWgood}) if $T^i=t$.
Note that in (a), because the coalescing random walks are stopped as soon as they meet
and are independent up until they meet, we can apply Corollary \ref{cor1}.
\er

We first consider the consequence of our definition of $T^i$: we define the random variables
$C_{i,j}$, $0\leq i<j\leq R$ by
\be{CijDef}
C_{i,j}=P\Big(X^{j}(v)=X^i(v) \text{ for some } v\in\big[t_j,T^i\big]
\,\,\big\vert\,\, X^i,\, T^i\Big).
\ee
We have for any $j\in\{i+1,i+2,\dots,R\}$
that
\be{a-cons5}
\begin{aligned}
C_{i,j}
&=
P\Big(X^{j}(v)=X^i(v) \text{ for some } v\in\big[t_i,T^i\big) \,\,\big\vert\,\,
X^i,\, T^i\Big)\\
&\quad +P\Big(X^{j}(T^i) = X^i(T^i),\, X^{j}(v) \neq X^i(v)\,\, \forall\, v<T^i
\,\,\big\vert\,\, X^i,\, T^i\Big).
\end{aligned}
\ee
By the definition of $T^{i,b}$, $|X^i(T^i)|\geq t^{1/3}-1$ unless $T^i< t_{i+1}$,
in which case $\{(X^{j}(s),s) \colon s\geq t_j\}$ cannot hit
$\{(X^i(u),u) \colon t_i\leq u\leq T^i\}$. Therefore, using a simple bound for
$p_s(\,\cdot\,)$ (use e.g.\ continuous version of Lawler \cite{law91}, Theorem 1.2.1,
inequality (1.10)), there exists some universal $K>0$ so that
\be{a-cons7}
\begin{aligned}
P\Big(X^{j}(T^i) = X^i(T^i),\, X^{j}(v) \neq X^i(v)\,\, \forall\, v<T^i
\,\,\big\vert\,\, X^i,\, T^i\Big)
&\leq \sup_{|x|\geq t^{1/3}-1}\,\,\sup_{u\geq 0} p_u(x)\\
&\leq \frac{K}{\big(t^{1/3}-1\big)^{2}}.
\end{aligned}
\ee
\br{rm4}
It is above all here we see the validity of the of the definition of $T^{i,b}$,
since this assures that for any $t_{i+1}\leq s\leq T^i$, $|X^i(s)|\geq t^{1/3}-1$.
Obviously the $1$ is arbitrary and could be replace by any $\lambda>0$.
\er
Combining (\ref{a-cons5}--\ref{a-cons7}) and summing over $i\leq j\leq R$ with $i<R$,
we obtain (recalling (a))
\be{CijProp}
\begin{aligned}
\sum_{j=i+1}^{R} C_{i,j}
&\leq \frac{\alpha}{10^2}+ \frac{RK}{(t^{1/3}-1)^2}\\
&<\frac{\alpha}{99},
\end{aligned}
\ee
for $t$ sufficiently large.

\bd{}
We say $1\leq j\leq R$ is \embf{good} if
\be{jGood}
\sum_{i=0}^{j-1} C_{i,j}\leq \frac{2\alpha}{99}.
\ee
\ed
\bl{GoodProplem}
At least $R/2$ of the $j$ are good.
\el
\bpr
By (\ref{CijProp}), we have
\be{GoodProp5}
\sum_{i=0}^{R-1}\sum_{j=i+1}^{R} C_{i,j}
\leq \frac{R\alpha}{99}.
\ee
Thus,
\be{GoodProp7}
\sum_{j=1}^{R}\sum_{i=0}^{j-1} C_{i,j}
\leq \frac{R\alpha}{99},
\ee
from which we obtain the result.
\epr
\bd{GoodRW}
We say a random walk $\big\{X^j(t_j +u)\colon u\geq 0\big\}$
is \embf{successful} if
\begin{itemize}
\item[(i)] the stopping time $T^j$ is equal to $t$;
\item[(ii)] $X^j$ does not hit a previous stopped random walk, i.e.,
for all $i<j$ and $s\in[t_j,T^i]$, $X^i(s)\neq X^j(s)$.
\end{itemize}
\ed

We consider now a somewhat unnatural filtration $\cF_0,\cF_1,\dots,\cF_R$.
Each of whose $\sigma$-fields will be based on the Poisson processes generating
the coalescing random walks. They are defined in the following way :
$\cF_0$ is trivial; $\cF_1$ is the $\sigma$-field generated by
$(X^0(u), 0\leq u\leq T^0)$; $\cF_r$ with $2\leq r\leq R$ is the $\sigma$-field
generated by $\cF_{r-1}$ and the random walk $X^{r-1}$ stopped at
$T^{r-1}\vee S^{r-1}$, where $S^{r-1}$ is the first time $(X^{r-1}(u),u)$
hits a previous (stopped) random walk. One way to see $\cF_r$ is as the $\sigma$-field
generated by the Harris system viewed along the paths of the $X^i$, $i\leq r-1$, that
is to say with information on $N^{x,y}$ for all $y$ on interval $I$ for $X^i(s)=x$
on $I$. It is clearly seen that on the $\sigma$-field $\cF_{j}$, the law of $(X^j(s),s)$
is simply a space-time random walk which evolves until it hits a point $(y,s)$ such that
$X^i(s)=y$ for some $i<j$ and $s\leq T^i$.

\bc{cor3}
If $t$ is sufficiently large, for at least $R/2$ random walks $X^j$, $1\leq j\leq R$
\be{GoodProba}
P\Big(X^j \text{ is successful }\,\,\big\vert\,\, \cF_j\Big)
\geq \frac{\alpha}{2}.
\ee
\ec
\bpr
By the definition of ``being good'' and Lemma \ref{GoodProplem}, for at least $R/2$
random walks $X^j$, $1\leq j\leq R$ we have $\sum_{i=0}^{j-1}C_{ij}\leq 2\alpha/99$.
Therefore, for those $j$,
\be{GoodProba5}
\begin{aligned}
P\Big(X^j \text{ hits a previous stopped random walk }\,\,\big\vert\,\, \cF_j\Big)
&\leq \sum_{i=0}^{j-1} P\Big(X^j\text{ hits }X^i\text{ stopped} \,\,\big\vert\,\, \cF_j\Big)\\
&=\sum_{i=0}^{j-1}C_{i,j}\leq \frac{2\alpha}{99}.
\end{aligned}
\ee
By Corollary \ref{cor2}, it follows that if $j$ is good
\be{GoodProba7}
P\Big(X^j \text{ is successful }\,\,\big\vert\,\, \cF_j\Big)
\geq \frac{2\alpha}{3}-\frac{2\alpha}{99}\geq \frac{\alpha}{2}>0.
\ee
\epr

As a consequence, we have the following result.
\bc{cor4}
There exists $K_6>0$ not depending on $t$ so that
\be{A1bound}
P\Big(\text{at least }K_6\log(t)\text{ random walks }\big(X^j \colon 1\leq j\leq
R\big)\text{ are successful}\Big)\geq 1- e^{-K_6\log(t)}.
\ee
In consequence, for the system $\chi^t$, except for an event of probability at most
$\exp[-K_6\log(t)]$, there exist at least $3t/4\leq s_1<s_2<\cdots<s_{\lfloor K_6\log(t)\rfloor}
\leq t$, such that
\begin{itemize}
\item[(i)] $\chi_u^{s_j}\neq\chi_{s_k-s_j+u}^{s_k}$ for all $1\leq j<k\leq K_6\log(t)$
and $0\leq u\leq s_j$;
\item[(ii)] $|\chi_u^{s_j}|\in (\sqrt{t},\sqrt{2t})$ for all $1\leq j\leq K_6\log(t)$
and $s_j-t/2 \leq u\leq s_j$.
\end{itemize}
\ec
\bpr
By Corollary \ref{cor3}, at least $R/2$ of the $1\leq j\leq R$ satisfy (\ref{GoodProba}).
For notational convenience only, we assume that (\ref{GoodProba}) holds for $1\leq j\leq R/2$.
Let $Z_j=\one{\{X^j\text{ is successful}\}}$. Therefore,
\be{A1-5}
P\big(Z_j=1 \,\vert\, Z_1, Z_2,\cdots,Z_{j-1}\big)
\geq \frac{\alpha}{2}\quad \forall\,1\leq j\leq \frac{R}{2}.
\ee
It follows that
\be{A1-7}
P\Big(\text{at least }\alpha R/8 \text{ random walks } \big(X^j \colon 1\leq j\leq R\big)
\text{ are successful}\Big)
\geq P\left(\sum_{i=1}^{R/2}Z_i\geq\frac{\alpha R}{8}\right).
\ee
We suppose that $(U_j\colon 1\leq j\leq R/2)$ is an i.i.d.\ sequence with uniform distribution
$\mathcal{U}([0,1])$ such that independently of the Harris system
\be{Ydef}
Y_j=Z_j\, \one\{U_j\leq \alpha/(2P(Z_j=1 \,\vert\, Z_1, Z_2,\cdots,Z_{j-1}))\}.
\ee
Therefore, $\big(Y_j \colon 1\leq j\leq R/2\big)$ a sequence of i.i.d.\ random variables
on $\{0,1\}$ so that
\be{A1-9}
P\big(Y_j=1\big)=\frac{\alpha}{2}
\quad\text{and}\quad
Y_j\leq Z_j\quad \forall\, 1\leq j\leq R/2.
\ee
Therefore,
\be{A1-11}
P\left(\sum_{i=1}^{R/2}Z_i\geq\frac{\alpha R}{8}\right)
\geq P\left(\sum_{i=1}^{R/2}Y_i\geq\frac{\alpha R}{8}\right).
\ee
But, by large deviations bound for Binomial process (see e.g. den Hollander \cite{hol00},
Chapter 1) and (\ref{Rdef}), we have
\be{A1-13}
\begin{aligned}
P\left(\sum_{i=1}^{R/2}Y_i\geq\frac{\alpha R}{8}\right)
&\geq 1-e^{-K\frac{\alpha}{4}R}\\
&\geq 1-e^{-\frac{K\alpha}{20C}\log(t)}
\end{aligned}
\ee
for some universal $K>0$ (not depending on $t$). Combining (\ref{A1-7}) and
(\ref{A1-11}--\ref{A1-13}), and reducing constants if necessary, we arrive at
(\ref{A1bound}).
\epr

\bprp{Prop1}
Let $K_1$ be a small positive constant to be more fully specified later.
Consider for all $0\leq i\leq K_1\log(t)$ the events $A_i(t)=$
\be{prop1-3}
\begin{array}{rll}
\Big\{&\hbox{there exist at least
       $3\times 2^{-i-2}t\leq s_1<s_2<\cdots<s_{\lfloor K_1\log(2^{-i}t)\rfloor}\leq 2^{-i}t$,
       such that}&\\
      &\hbox{(i) $\chi_u^{s_j}\neq\chi_{s_k-s_j+u}^{s_k}$ for all $1\leq j<k\leq K_1\log(2^{-i}t)$
       and $0\leq u\leq s_j$;}&\\
      &\hbox{(ii) $|\chi_u^{s_j}|\in (\sqrt{2^{-i}t},\sqrt{2^{-i+1}t})$ for all
       $1\leq j\leq K_1\log(2^{-i}t)$ and $s_j-2^{-i-1}t\leq u\leq s_j$}
&\Big\}.
\end{array}
\ee
Thus, under this definition, Corollary \ref{cor4} says that
\be{prop1-5}
P\big(A_i(t)\big)\geq 1-\exp\Big[-K_1\log\big(2^{-i}t\big)\Big]
\geq 1-\exp\bigg[-\frac{K_1}{2}\log(t)\bigg]
\ee
if $K_1$ is small enough.
Therefore, we have that (after reducing $K_1$)
\begin{itemize}
\item[(i)] events $A_i(t)$ are independent for $1\leq i\leq K_1\log(t)$;
\item[(ii)] $P\big(A_i(t)\big)\geq 1-\exp\big[-K_1\log(t)\big]$.
\end{itemize}
If $\sum_{i\leq K_1 \log(t)}\one_{A_i^{\c}}<K_1\log(t/2)$, then $\#\chi^t \geq K_1\log^2(t)$.
Therefore, there exists $K_2>0$ so that
\be{prop1-7}
\begin{aligned}
P\Big(\#\chi^t\leq K_1\log^2(t)\Big)
&\leq P\left(\sum_{i\leq K_1\log(t)}\one_{A_i^{\c}}\geq K_1\log\bigg(\frac{t}{2}\bigg)\right)\\
&\leq 2^{K_1\log(t)}\, \exp\bigg[-\frac{K_1^2}{2}\log^2(t)\bigg]
\leq e^{-K_2 \log^2(t)}.
\end{aligned}
\ee
\eprp

%%%%%%%%%%%%% SUBSECTION 3.3 %%%%%%%%%%%%%%%%%%%%%%%%%%%%%%%%%%%%%%%%%%%%%%%%%%%%%%%%%%%%%%%%%%%%%%

\subsection{Proof of Theorem \ref{Th2}}
\label{S3.3}

Denote by $\#\chi^{[r,s]}$, $0\leq r\leq s$, the number of distinct sites in the collection
$\{\chi_u^u \colon r\leq u\leq s\}$. We refer to Bramson, Cox and Griffeath \cite{bramcoxgri88},
Section 2:
\bl{BCGlem}
There exists some positive finite constant $K$ so that for all $t>1$
\be{BCG1}
E\Big(\#\, \chi^{[t/2,t]}\Big)\leq K\log(t).
\ee
\el
We are now ready to prove Theorem \ref{Th2}.

\bpr
For all $\alpha\in(\rho,1)$ and $t\geq 0$, by Jensen's inequality, we have
\be{lb3}
\begin{aligned}
\P_{\nu_\rho}\big(T_t\geq \alpha t\big)
&\geq \P_{\,\nu_\rho}\bigg(\int_{(1-\alpha)t}^{t}\one\{\eta_s(0)=1\}\, ds
=\alpha t\bigg)\\
&\geq \rho^{E\big(\# \chi^{[(1-\alpha)t,t]}\big)}.
\end{aligned}
\ee
Split time interval $((1-\alpha)t,t]$ so that
\be{lb5}
((1-\alpha)t,t]
\subset\bigcup_{k=0}^{\lfloor-\log_2(1-\alpha)\rfloor}
\big(t2^{-k-1},t2^{-k}\big],
\ee
then apply Lemma \ref{BCGlem} to each $\chi^{[t2^{-k-1},t2^{-k}]}$,
$k=0,\ldots,\lfloor-\log_2(1-\alpha)\rfloor$ to obtain
\be{lb7}
E\Big(\#\, \chi^{[(1-\alpha)t,t]}\Big)
\leq K_1\log(t)
\ee
for $K_1$ a finite positive constant large enough.
Then, combining (\ref{lb3}) and (\ref{lb7}), we get
\be{lb11}
\P_{\nu_\rho}\big(T_t\geq \alpha t\big)
\geq e^{-C_1\log(t)}
\ee
for $C_1$ large enough.
\epr
%%%%%%%%%%%% REFERENCES %%%%%%%%%%%%%%%%%%%%%%%%%%%%%%%%%%%%%%%%%%%%%%%%%%%%%%%%%%%%%%%%%%%%%%%%%%%

\end{document}